\documentclass[11pt]{article}

\usepackage{amsfonts}
\usepackage{amsthm}
\usepackage{amsmath}
\usepackage{amsfonts}
\usepackage{latexsym}
\usepackage{amssymb}
\usepackage[dvips]{color}

 \newtheorem{thm}{Theorem}[subsection]
 \newtheorem{cor}[thm]{Corollary}
 \newtheorem{lem}[thm]{Lemma}
 \newtheorem{prop}[thm]{Proposition}
 \newtheorem*{conj}{Conjecture}
 \newtheorem{defn}[thm]{Definition}
 \newtheorem{rem}[thm]{Remark}
 
 \newtheorem{exa}[thm]{Example} 

 \numberwithin{equation}{section}

\def\noi{\noindent}
\definecolor{color}{RGB}{140,140,140}
\def\QED{\hfill $\color{color}{\blacksquare}$}
\def\EOE{\hfill $\color{color}{\blacktriangle}$}
\def\bdem{\proof}
\def\edem{\renewcommand{\qed}{\hfill $\color{color}{\blacksquare}$}}

\def\eps{\varepsilon}

\def\la{\lambda}
\def\cF{\mathcal{F}}
\def\ga{\gamma}

\def\ds{\displaystyle} 

\def\rai{^{1/2}}
\def\mrai{^{-1/2}}

\def\N{\mathbb{N}}

\def\R{\mathbb{R}}
\def\C{\mathbb{C}}
\def\inc{\subseteq}
\def\bm{\left(\begin{array}}
\def\em{\end{array}\right)}

\def\orto{^\perp}

\def\cB{\mathcal{B}}

\def\cD{\mathcal{D}}
\def\cE{\mathcal{E}}

\def\cM{\mathcal{M}}

\def\cU{\mathcal{U}}
\def\cV{\mathcal{V}}
\def\cW{\mathcal{W}}

\def\ese{\mathcal{S}}

\def\ewe{\mathcal{W}}

\def\ben{\begin{enumerate}}
\def\een{\end{enumerate}}
\def\beq{\begin{equation}}
\def\eeq{\end{equation}}
\def\barr{\begin{array}}
\def\earr{\end{array}}
\def\inv{^{-1}}
\def\glh{\mathcal{G}\textit{l} \, (\hil)}

\newcommand{\peso}[1]{ \quad \text{ #1 } \quad }
\newcommand{\sub}[2]{{#1}_{_{#2}}}

\DeclareMathOperator{\Preal}{\R\mbox{e}} 

\DeclareMathOperator*{\dist}{dist}

\DeclareMathOperator{\tr}{tr}

\DeclareMathOperator{\sgn}{sgn}

\DeclareMathOperator{\leqp}{\leqslant}

\def\ga{\gamma}

\newcommand{\diag}[1]{\hbox{\rm diag}\left( #1\right)}

\newcommand{\pint}[1]{\displaystyle \left \langle #1 \right\rangle}

\newcommand{\hil}{\mathcal{H}}
\newcommand{\op}{L(\mathcal{H})}
\newcommand{\opsa}{L_{sa}(\mathcal{H})}
\newcommand{\posop}{L(\mathcal{H})^+}

\newcommand{\cene}{\mathbb{C}^r}
\newcommand{\mat}{\mathcal{M}_r (\C) }
\newcommand{\matu}{\mathcal{U}(r)}
\newcommand{\matsa}{\mathcal{M}_r^{h}(\C)  }
\newcommand{\matah}{\mathcal{M}_r^{ah}(\C)  }

\newcommand{\matinv}{\mathcal{G}\textit{l}\,_r(\C) }

\newcommand{\spec}[1]{\sigma\left( #1\right)}

\newcommand{\conv}{\xrightarrow[n\rightarrow\infty]{}}

\newcommand{\alus}[1]{\Delta_{1/2}\left(#1\right)}
\newcommand{\alusit}[2]{\Delta_{1/2}^{#1}\left(#2\right)}

\newcommand{\alulf}[1]{\left|#1\right|^{\la}U\left|#1\right|^{1-\la}}
\newcommand{\aluit}[2]{\Delta^{#1}\left(#2\right)}

\newcommand{\alul}[1]{\Delta_\la\left(#1\right)}
\newcommand{\alum}[1]{\Delta_{1/2}\left(#1\right)}

\newcommand{\alulit}[2]{\Delta_\la^{#1}\left(#2\right)}

\def\iiff{ if and only if }

\newcommand{\orb}[1]{\ese \left({#1}\right)}
\newcommand{\orbu}[1]{\cU \left({#1}\right)}
\newcommand{\orbl}[1]{\ese_{L} \left({#1}\right)}
\newcommand{\orbul}[1]{\cU_{L} \left({#1}\right)}
\newcommand{\orbla}[1]{\ese_{o} \left({#1}\right)}
\newcommand{\orbula}[1]{\cU_{o} \left({#1}\right)}
\newcommand{\der}[3]{\sub{T}{#2}{#1} \left({#3}\right)}
\newcommand{\dersin}[2]{\sub{T}{#2}{#1}}
\newcommand{\had}[1]{\Psi_{#1}}
\newcommand{\preal}{\sub{P}{\R e}}
\newcommand{\pim}{\sub{P}{\mathbb{I}{m}}}

\newcommand{\qd}{\sub{Q}{D}}
\newcommand{\qn}{\sub{Q}{N}}
\newcommand{\diagr}{\cD_{r}(\C)}
\newcommand{\diagdis}{\cD_{r}^*(\C)}

\newcommand{\rk}[1]{\mbox{rk #1 }}

\begin{document}

\title{
Convergence of iterated Aluthge transform \\ sequence for diagonalizable matrices II: \\ 
 $\la$-Aluthge transform 
\author{Jorge Antezana
\thanks{\noi $^*$Partially supported by CONICET (PIP 4463/96), Universidad de La Plata (UNLP 11 X472) and
ANPCYT (PICT03-09521).}
\and Enrique R. Pujals 
\thanks {\noi Partially supported by CNPq} 
\and Demetrio Stojanoff $^*$ 
}
\date{}
}
\maketitle

\noi {\bf Jorge Antezana and Demetrio Stojanoff}

\noi Depto. de Matem\'atica, FCE-UNLP,  La Plata, Argentina
and IAM-CONICET  

\noi e-mail: antezana@mate.unlp.edu.ar and demetrio@mate.unlp.edu.ar

\medskip

\noi {\bf Enrique R. Pujals}

\noi Instituto Nacional de Matem\'atica Pura y Aplicada (IMPA), Rio de Janeiro, Brasil.

\noi e-mail: enrique@impa.br

\noi

\medskip
\noi

\vglue1truecm

\begin{abstract}
Let $\lambda \in (0,1)$ and let $T$ be a $r\times r$ complex matrix with polar decomposition $T=U|T|$. Then, the $\la$- Aluthge transform is defined by 
$$
\Delta_\lambda \left(T \right)= |T|^{\lambda} U |T |^{1-\lambda}.
$$ 
Let $\Delta_\lambda^{n}(T)$ denote the n-times iterated Aluthge transform of $T$, $n\in\mathbb{N}$. We prove that the sequence $\{\Delta_\lambda^{n}(T)\}_{n\in\mathbb{N}}$ converges for every $r\times r$ {\bf diagonalizable}  matrix $T$.
We show regularity results for the two parameter map $(\la , T) \mapsto \alulit{\infty}{T}$, and we study for which matrices the map $(0,1)\ni \lambda  \mapsto \Delta_\lambda^{\infty}( T)$ is constant.
\end{abstract}

\newpage

\section{Introduction.}
Let $\hil$ be a complex Hilbert space, and let $\op$ be the
algebra of bounded linear operators on $\hil$. Given $T\in\op$,
consider its (left) polar decomposition $T=U|T|$. In order to
study the relationship among  p-hyponormal operators, Aluthge
introduced in \cite {[Aluthge]} the transformation
$\alum{\cdot}:\op\to\op$ defined by
\[
\alus{T}  = |T|\rai U |T|\rai.
\]

\noi Later on, this transformation, now called Aluthge transform, 
was also studied in other contexts by several authors, such
as Jung, Ko and Pearcy \cite{[JKP0]} and  \cite{[JKP1]}, 
Ando \cite{[Ando]}, Ando and Yamazaki \cite
{[Ando-Yamaza]}, Yamazaki \cite{[Yamazaki]}, Okubo
\cite{[lambdaluthges]} 
and Wu \cite {[Wu]} among others.

\medskip
\noi In this paper, given $\la\in (0,1)$ and $T\in\op$, we study
the so-called $\la$-Aluthge transform of $T$ defined by
\[
\alul{T}=\alulf{T}.
\]
This notion has already been considered by Okubo in
\cite{[lambdaluthges]} (see also \cite{[AMS]} and \cite {[HT]}). 
We denote by  $\alulit{n}{T}$ the n-times iterated $\la$-Aluthge transform of $T$, i.e.
\begin{align}\label{laal}
\alulit{0}{T}=T; \peso{and}
\alulit{n}{T}=\alul{\alulit{n-1}{T}}\quad n\in\N.
\end{align}
In a previous paper \cite{[APS]}, we show that the iterates of usual Aluthge transform 
$\alusit{n}{T}$ converge to a normal matrix $\alusit{\infty}{T}$ for every diagonalizable 
matrix $T \in \mat$ (of any size). We also proved in \cite{[APS]} the smoothness of the map 
$T \mapsto \alusit{\infty}{T}$ when it is restricted to a similarity orbit, or to the 
(open and dense) set $\diagdis$ of invertible $r\times r$ matrices  with $r$ different 
eigenvalues. The key idea  was  to use  a dynamical systems approach
to the Aluthge transform, thought as acting on the similarity orbit 
of a diagonal invertible matrix. 
Recently, Huajun Huang and Tin-Yau Tam \cite{[HT]} showed, with other approach, 
that the iterates of every $\la$-Aluthge transform 
$\alulit{n}{T}$ converge,  for every matrix $T \in \mat$ with all its eigenvalues of 
{\bf different moduli}.

In this paper, we study the general case of $\la$-Aluthge transforms 
by means of a dynamical systems approach. 
This allows us to generalize 
Huajun Huang and Tin-Yau Tam result 
for every {\bf diagonalizable} matrix $T \in \mat$,  as well as to
show  regularity results for the 
two parameter map $(\la , T) \mapsto \alulit{\infty}{T}= \lim_{n\in \N} \alulit{n}{T}$. 

Now we briefly describe 
the dynamical point of view of the problem:
For every $\la \in (0,1)$ and any invertible matrix $T$, it holds that 
$
\alul{T}=|T|^{\la } \, T\, |T|^{-\la}.
$  
So the $\la$-Aluthge transform of $T$ belongs to the similarity orbit of $T$. This suggests  
that we can study the Aluthge transform restricted to the similarity orbit of some invertible operator. 
From that point of view, the diagonalizable case has a better dynamical behavior, as 
detailed in \cite{[APS]}. 
If $T$ is diagonalizable, the similarity orbit of $T$ coincides with the similarity orbit of some 
diagonal operator $D$, which we denote $\orb{D}$. The unitary orbit of $D$, denoted by $\orbu{D}$, 
is a compact submanifold of $\orb{D}$ that consists of all normal matrices in $\orb{D}$. Hence 
$\orbu{D}$ is  fixed by the Aluthge transform and, as it was shown in \cite{[AMS]}, 
all the limit points of the sequence $\{\alulit{n}{T}\}_{n\in\N}$ belong to $\orbu{D}$. 

As  it was shown in \cite {[APS]} for $\la=1/2$, we 
show that for any $N\in \orbu{D}$ there is a local submanifold 
$\ewe_{\la, N}^s$ transversal to $\orbu{D}$ characterized as the set the matrices 
(near $N$) that converges with a exponential rate to $N$ by the iteration of the 
$\la$-Aluthge transform. Moreover, the union of these submanifolds form an open 
neighborhood of $\orbu{D}$ (see Corollary \ref{entorno y algo mas}). 
Thus, since the sequence $\{\alulit{n}{T}\}_{n\in\N}$ goes toward $\orbu{D}$, 
for some $n_0$ large enough the sequence of iterated Aluthge transforms gets into
this open neighborhood and converge exponentially. 

These results follow from the classical arguments of stable manifolds (first introduced independently  
by Hadamard and Perron, see Theorem \ref{teorema 5.5}; for details and 
general results 
about the stable manifold theorem, see \cite{[HPS]} or the Appendix of \cite{[APS]}). 
In order to use the stable manifold theorem, we show previously that the 
derivative of the  $\la$-Aluthge transform in any $N\in\orbu{D}$ has two invariant complementary 
directions, one tangent to $\orbu{D}$, and other one transversal to it where the derivative 
is a contraction (see Theorem \ref{EL TEO}). Using these techniques  
we prove that $\alulit{n}{T} \conv \alulit{\infty}{T} \in \matu$, 
for every $r\times r$  diagonalizable  matrix $T$. 
 We also prove that the two parameter map 
 $\Delta_{\infty} (\la , T) = \alulit{\infty}{T}$ is of class $C^\infty$, when restricted to
$(0,1) \times \diagdis$ and to $(0,1) \times \orb{D}$ for a diagonal matrix $D \in \mat$. 

We also study for which matrices $T$ the map 
$(0,1)\ni \la  \mapsto R_T(\la) = \alulit{\infty}{T}$ is constant. Some partial results are obtained, 
in particular that $R_T$ is not constant for most diagonalizable matrices $T$. We also show that 
$R_T$ is constant for every $T$ in the similarity 
orbit of a diagonal matrix $D$ if  $\sigma (D) = \{ d_1 \, , d_2\}$ with  $|d_1 | = |d_2 |$. 
We state and discuss the following conjecture: this is the unique case where 
$R_T$ is constant for every $T$ in the similarity 
orbit of  $D$.

The paper has a structure very similar to \cite{[APS]} because, at any step of the 
dynamical systems argument, 
we need to state results which differ slightly of those results of \cite{[APS]}.  
The  proofs of these results are omitted or just sketched. 
The paper is organized as follows: 
in section 2, we collect several preliminary definitions and results about the
the stable manifold theorem, about the geometry of similarity and unitary orbits, 
and about known results on $\la$-Aluthge transforms. 
In section 3, we prove the convergence results. In section 4 we study the smoothness of 
the two parameter map $(\la , T) \mapsto \alulit{\infty}{T}$ and we study the behavior  
of the limit function $\alulit{\infty}{T}$ with respect  to the parameter $\la$. 
The basic tool for these results, in order to apply the stable manifold theorem to the 
similarity orbit of a diagonal matrix, is the mentioned  
Theorem \ref{EL TEO}, whose proof, somewhat technical, is done in 
section 5. 

\section{Preliminaries.}
In this paper $\mat$ denotes the algebra of complex $r\times r$ matrices, 
$\matinv$ the group of all invertible elements of $\mat$, $\matu$ the group 
of unitary operators, and $\matsa$ (resp. $\matah\, $) denotes the real algebra 
of Hermitian (resp. anti Hermitian) matrices. 
Given $T \in \mat$, $R(T)$ denotes the
range or image of $T$, $\ker T $ the null space of $T$, 
$\sigma (T)$ the spectrum of $T$, $\rk{T} = \dim R(T)$ the rank of $T$, 
$\tr(T)$ the trace of $T$,
and $T^*$ the adjoint of $T$. If $v \in \cene$, we denote by 
$\mbox{\rm diag}(v) \in \mat$ the diagonal matrix with $v$ in its diagonal. 
We shall consider the space of matrices $\mat$ as a real Hilbert space 
with the inner product defined by
$$
\pint{A,\ B}=\Preal\big(\tr(B^*A)\big).
$$
The norm induced by this inner product is the so-called Frobenius norm, 
denoted by $\|\cdot \|_{_2}\,$. 

\medskip
\noi
On the other hand, let $M$ be a manifold. By means of $TM$ we denote the tangent bundle of $M$ and by means of $T_xM$ we denote the tangent space at the point $x\in M$. Given a function $f\in C^{r}(M)$, where $r=1,\ldots,\infty$, $\der{f}{x}{v}$ denotes the derivative of $f$ at the point $x$ applied to the vector $v$.

\subsection{Stable manifold theorem}

In this section we state the stable manifold theorem for an invariant  set of a smooth 
endomorphism (see \ref{teorema 5.5} below). We refer to \cite{[APS]} for a detailed 
description of these type of results. 
Let $M$ be a smooth Riemann manifold and $N\subseteq M$ a submanifold (not 
necessarily compact). Throughout this subsection $\sub{T}{N}M$ denotes the 
tangent bundle of $M$ restricted to $N$.

\begin{defn}\label{prelamination}\rm
A $C^r$ \textit{pre-lamination} indexed by $N$ is a continuous choice of a $C^r$ embedded disc $\cB_x$ 
through each $x\in N$. Continuity means that $N$ is covered by open sets  $\cU$ 
in which $x\to B_x$ is given by
$$
\cB_x=\sigma(x)((-\eps,\eps)^k)
$$
where $\sigma: \cU \cap N\to \mbox{Emb}^r((-\eps,\eps)^k,M)$ is a continuous section. 
Note that 
$\mbox{Emb}^r((-\eps,\eps)^k,M)$ is a $C^r$ fiber bundle over $M$ whose 
projection is $\beta\to \beta(0)$. Thus $\sigma(x)(0)=x$. 
If the sections mentioned above are $C^s$, $1\leq s\leq r$, we say that the $C^r$ 
pre-lamination is of class $C^s$. \EOE
\end{defn}

\begin{defn}\label{pseudohyperbolic}\rm
Let $f$ be a smooth endomorphism of $M$, $\rho>0$, and suppose that $\left. f\right|_{N}$ is a homeomorphism. Then, $N$ is  \textit{$\rho$-pseudo hyperbolic} for $f$ if there exist two smooth subbundles of $\sub{T}{N}M$, denoted by $\cE^s$ and $\cF$, such that 
\begin{enumerate}
	\item $\sub{T}{N}M= \cE^s \oplus \cF$;
        \item $\sub{T}{}N=  \cF$;
	\item Both, $\cE^s$ and $\cF$, are $Tf$-invariant;
	\item $T\, f$ restricted to $\cF$ is an automorphism, which expand it by a factor
	greater than $\rho$. 
	\item $\dersin{f}{x}:\cE_x^s\to\cE_{f(x)}^s$ has norm lower than $\rho$. \EOE
\end{enumerate}
\end{defn}

\noi
In this case, the stable manifold theorem assures that for any point $x\in N$ it 
is possible to find an $f-$invariant submanifold transversal to $N$ tangent to $\cE^s$ 
and characterized as the set of points with trajectories asymptotic to the trajectory of $x$. 
We shall state the following particular version of this theorem. 
For a proof, see Theorem 2.1.4, Corollary 2.1.5 
and Remark 2.1.6 of \cite{[APS]}.

\begin{thm}[Stable manifold theorem for a submanifold of fixed points]\label{teorema 5.5} \rm 
Let $f$ be a $C^r$ endomorphism of $M$ with a $\rho$-pseudo hyperbolic submanifold 
$N$ with  $\rho< 1$. Assume that any point $p$ in $N$ is a {\bf fixed point}.
Then, there is a $f$-invariant $C^r$-pre-lamination  
$\ewe^s:  N\to \mbox{Emb}^r((-1,1)^k,M)$ {\bf of class $\mathbf {C^r}$}  
such that, for every $x \in N$, 
\begin{enumerate}
\item $\ewe^s(x)(0)=x$,
\item  $\ewe_x^s=\ewe^s(x)((-1,1)^k)$  is tangent to $\cE_x^s$  at every $x\in N$,
\item There exists an open neighborhood $\cU$ 
of $x$ (open relative to $M$) such that
\begin{equation}\label{converge}
	\ewe_x^s\cap \cU=\Big\{y\in \cU:\ \dist(x,f^n(y))< \  \dist(x,y)\,  \rho^n \Big\}.
	\end{equation}
\item 	
	 If $y \in N$, $x \neq y$, then $\ewe_x^s\cap \ewe_y^s = \varnothing  $. 

\item There exists $\ga>0$ such that  
	$$B(x,\ga)\subset \bigcup_{x\in N} \ewe_x^s \ .	$$
This implies that $\displaystyle \bigcup_{x\in N} \ewe_x^s$ contains 
an open neighborhood $\ewe(N)$ of $N$ in $M$. 
	\item The map 
\begin{equation}\label{proyecta}
p : \ewe(N)  \to N \peso{ given by }
p(a) = x \peso{if} a \in  \ewe_x^s(x) \ , 
\end{equation}
is well defined  and it is of class $C^r$. \QED
\end{enumerate}

\end{thm}

\subsection{Similarity orbit of a diagonal matrix}

In this subsection we recall some facts about the similarity orbit of a diagonal matrix. 

\begin{defn}\rm
Let  $D\in\mat$ be diagonal. 
The similarity orbit of $D$ is the set   
$$
\orb{D} = \{ \ SDS\inv\  : \ S \in \matinv \ \} \ . 
$$
On the other hand, $\orbu{D}= \{ \ UDU^* \  : \ U \in \matu \ \}$ 
denotes the unitary orbit of $D$. 
We donote by 
$\sub{\pi}{D} : \matinv \to \orb{D} \inc \mat $  the $C^\infty$ map defined by 
$\sub{\pi}{D}(S) = SDS\inv $.  
With the same name we note its restriction to the unitary group:  
$\sub{\pi}{D} : \matu \to \orbu{D} $. 
\EOE
\end{defn}

\begin{prop}\label{son variedades} \rm
The similarity orbit $\orb{D}$ is a $C^\infty$ submanifold of $\mat$, and the projection
$\sub{\pi}{D} : \matinv \to \orb{D}$ becomes a submersion. Moreover, $\orbu{D}$ is a compact submanifold of $\orb{D}$, which consists of the normal elements of $\orb{D}$, and 
$\sub{\pi}{D} : \matu \to \orbu{D}$ is a submersion. In particular, the maps 
$\sub{\pi}{D}$ have $C^\infty$ (similarity and unitary) local cross sections.  \QED
\end{prop}

\noi
For every $N=UDU^*\in\orbu{D}$, it is well known (and easy to see) that
\begin{align}
\sub{T}{N}\,\orb{D}&= \sub{T}{I}( \pi_N ) (\mat\,)  =\{[A,N]=AN-NA: \ A\in\mat\}. \nonumber\\
\intertext{In particular}
\sub{T}{D}\,\orb{D}&=\{AD-DA: \ A\in\mat\} \nonumber
\\\label{tan1} &=\{X\in\mat:\ X_{ij}=0 \ \mbox{for every $(i,j)$ such that $d_i=d_j$}\}.\\
\intertext{Note that, }
\sub{T}{N}\,\orb{D}&=\{[A,N]=AN-NA: \ A\in\mat\}\nonumber
\\&=\{(UBU^*)UDU^*-UDU^*(UBU^*): \ B\in\mat\}\nonumber
\\&=\{U[B,D]U^*=BD-DB: \ B\in\mat\}=U\Big(\sub{T}{D}\, \orb{D}\Big)U^*\ .\\
\intertext{On the other hand, since $\sub{T}{I}\, \matu = \matah= \{ A \in \mat : A^* = -A\}\,$, we obtain}
\sub{T}{D}\,\orbu{D}
&= \sub{T}{I}( \pi_D ) (\matah \,) =\{[A,D]=AD-DA: \ A\in\matah\} \peso{and} 
\ ,\nonumber\\
\label{ss}\sub{T}{N}\,\orbu{D}&=\{[A,N]=AN-NA: \ A\in\matah\}=U\Big(\sub{T}{D}\,\orbu{D}\Big) U^* \ . 
\end{align}

\noi Finally, along this paper we shall consider on $\orb{D}$ (and in $\orbu{D}$) the 
Riemannian structure inherited from $\mat$ (using the usual inner product on their 
tangent spaces). For $S, T \in \orb{D}$, we denote by $\dist(S, T)$  the Riemannian 
distance between $S$ and $T$ (in $\orb{D}\,$). 
Observe that, for every $U \in \matu$, one has that $U\orb{D} U^* = \orb{D} $ and the map 
$T \mapsto UTU^*$ is isometric, on $\orb{D}$,  with  respect to the Riemannian metric as 
well as with respect to the $\| \cdot \|_{_2}$ metric of $\mat$.

\subsection{$\la$-Aluthge transforms}

\begin{defn}\rm
Let $T\in\mat$, and suppose that $T=U|T|$ is the polar
decomposition of $T$. Let $\la \in (0,1)$. Then, we define the
$\la$-Aluthge transform of $T$ in the following way:
\begin{align*}
\alul{T}&=\alulf{T}
\end{align*}
We denote by  $\alulit{n}{T}$ the n-times iterated $\la$-Aluthge transform of $T$, i.e.
\begin{align*}
\alulit{0}{T}=T; \peso{and}
\alulit{n}{T}=\alul{\alulit{n-1}{T}}\quad n\in\N.
\end{align*}
\end{defn}

\noi The following proposition contains some properties of $\la$-Aluthge transforms which 
follows easily from its definition.

\begin{prop}\label{facilongas} \rm
Let $T\in\mat$ and $\la \in (0,1)$. Then:
\begin{enumerate}
  \item $\alul{cT}=c\alul{T}$ for every $c\in \C$.
  \item $\alul{VTV^*}=V\alul{T}V^*$ for every $V\in\matu$.
  \item If \ $T=T_1\oplus T_2$ then    
  $\alul{T}=\alul{T_1}\oplus \alul{T_2}$.
  \item $\|\alul{T}\|_{_2}\leqp \|T\|_{_2}$.
  \item $T$ and $\alul{T}$ have the same  characteristic polynomial. 
  \item In particular,   $\spec{\alul{T}}=\spec{T}$.
\end{enumerate}
\end{prop}

\noi The following theorem states the regularity properties of $\la$-Aluthge transforms.

\begin{thm}\label{continuidad}
The $\la$-Aluthge transform is 
continuous in $\mat$. 
Moreover, the map $(0, 1) \times \matinv \ni (\la , T) \mapsto \alul{T} $  is of class $C^\infty$.
\end{thm}

\proof 
The continuity part was proved in \cite{[AMS]} (see also \cite{[Dykema]}). If 
$T \in \matinv$, then 
$$
\alul{T} = |T|^\la \ T \ |T|^{- \la} \ , \quad \la \in (0,1) \ .
$$
This clearly implies regularity, since the map 
$$
(0, 1) \times \matinv \ni (\la , T)\ \   \longmapsto  \ \ |T|^\la = (T^*T)^{\la / 2} = 
\exp ( \la/2 \, \log T^*T ) 
$$ 
is of class $C^\infty $ . \QED

\medskip
\noi The following result is proved  in \cite{[AMS]}:

\begin{prop}\label{puntos limites normales} \rm
Given $T\in\mat$ and $\la \in (0,1)$, the limit points of the sequence
$\{\alulit{n}{T}\}_{n\in \N}$ are normal. Moreover, if $L$ is a limit point,
then $\spec{L}=\spec{T}$ with the same algebraic multiplicity. 
In particular, for each $\la \in (0,1)$, 
one has that  $\alul{T} = T $ if and only if $T$ is normal. 
\end{prop}

\noi Finally, we mention a result concerning the Jordan structure of 
Aluthge transforms proved in \cite{[AMS]}. We need the following definitions. 

\begin{defn}\rm Let $T \in \mat $ and $\mu \in \C$. We denote
\ben
\item $m(T, \mu)$ the $algebraic$ $multiplicity$ of $\mu$ for $T$, i.e. the maximum $n \in \N$ such that
$(x-\mu)^n $ divides the characteristic polynomial of $T$.  
\item $m_0(T, \mu)= \dim \ker( T-\mu I )$,
the $geometric$ $multiplicity$ of $\mu$. \EOE
\een
\end{defn}

\medskip
\noi
Observe that Proposition \ref{facilongas} says that 
$m (T, \mu ) =  m ( \alul{T}, \mu )$ for every $\mu \in \C$.  

\begin{prop}\label{algeo}
Let $T\in \mat$. 
\ben
\item If $\ 0 \in \spec{T}$, then, there exists $n\in\N$ such that
$$
 m(T, 0)=  m_0(\alulit{n}{T}, 0) = \dim \ker(  \alulit{n}{T}).
$$ 
\item For every  $\mu \in \sigma(T)$, 
$m_0(T, \mu ) \leqp m_0( \alul{T}, \mu ).$  
\een
\end{prop}

\noi
Observe that this implies that, if $T$ is diagonalizable (i.e. $m_0(T, \mu )  = m(T, \mu ) $ for every $\mu$), 
then also $\alul{T}$ is diagonalizable. 

\begin{rem}\rm
Let $T \in \matinv$ with polar decomposition $T = U|T|$. The Duggal or $1$-Aluthge transform of $T$ 
is $\Delta_1(T) = |T|U$. It is easy to see that the map $\Delta_1 : \matinv \to \matinv$ 
is continuous and that $\Delta_1(T) = T$ \iiff $T$ is normal. 
Observe that $U \in \matu$, so that  $\Delta_1(T) = UTU^*$, and the distance of $\Delta_1(T)$ 
to the normal matrices is the same as the 
distance of $T$ to the normal matrices. All these facts imply the iterated Duggal transforms 
$\Delta_1^n(T)$ can not converge, unless $T$ is normal. \EOE
\end{rem}

\section{Convergence}
In this section, we prove the convergence of iterated $\la$- Aluthge transforms for 
every diagonalizable matrix and $\la \in (0,1)$. Throughout the next subsections, a diagonal  matrix $D=\mbox{diag}(d_1,\ldots,d_n)\in\matinv$ is fixed.
For every $j\in\{1,\ldots,n\}$, let $d_j=e^{\,i\theta_j}|d_j|$ be the polar 
decomposition of $d_j$, where $\theta_j\in [0,2\pi]$.

\subsection{Reduction to the invertible case}
We claim that the proof of the convergence of iterated $\la$- Aluthge transforms can be reduced to 
the invertible case. Indeed, let  $T\in \mat$ be a a diagonalizable matrix with polar decomposition
$T=U|T|$. As $R(T)$ is a (oblique) complement of $\ker T= \ker |T|^\la$ and $R(U|T|^{1-\la})=R(T)$,  it holds that
$$
R(\alul{T})=R(\alulf{T})=R(|T|).
$$
On the other hand, it is easy to see that $\ker \alul{T}=\ker |T|^{1-\la}=\ker |T|$, which is orthogonal to $R(|T|)$. By Proposition \ref{algeo}, after one iteration we get that 
\beq\label{tuno}
\alul{T}=\begin{pmatrix}
T_1&0\\
0  &0
\end{pmatrix}\begin{array}{l}
\ker T ^\bot\\
\ker T
\end{array}\, ,
\end{equation}
where  $T_1$ is invertible and diagonalizable on $\ker T^\bot$.  
By Proposition \ref{algeo} again, 
$$
\alulit{n}{T}=\begin{pmatrix}
\alulit{n-1}{T_1}&0\\
0  &0
\end{pmatrix}\begin{array}{l}
\ker T ^\bot\\
\ker T 
\end{array}\ , \peso{for every $n\in \N $ .}
$$
Hence,  the convergence of $\{\alulit{n}{T}\}_{n\in\N}$ is equivalent to the convergence of $\{\alulit{n}{T_1}\}_{n\in\N}\,$.

\subsection{Main Theorem}
Reduced the problem to the invertible case, the key tool, which allows to use 
the stable manifold theorem \ref{teorema 5.5}, is Theorem \ref{EL TEO} below. The 
proof of this theorem is rather long and technical, 
for this reason, we postpone it until section \ref{la prueba}, and we 
continue in this section with its consequences. 

\begin{thm}\label{EL TEO} \rm
Let $\la \in (0,1)$. 
The $\la$-Aluthge transform $\alul{\cdot}:\orb{D}\to\orb{D}$ is a $C^\infty$ map, and for every $N\in\orbu{D}$, there exists a subspace $\sub{\cE}{N,\la}^s$ in the tangent space $\sub{T}{N}\orb{D}$ such that
\begin{enumerate}
	\item $\sub{T}{N}\orb{D}=\sub{\cE}{N,\la}^s\oplus \sub{T}{N}\orbu{D}$;
	\item Both, $\sub{\cE}{N,\la}^s$ and  $\sub{T}{N}\orbu{D}$, are $T_N\,\Delta_\la$-invariant;
	\item $\left\|\left. T_N\,\Delta_\la\right|_{\sub{\cE}{N,\la}^s}\right\|\leq\sub{k}{D,\la}<1$, where 
$$
\sub{k}{D,\la}  = \max\left\{ \max_{|d_i|\neq |d_j|}\frac{|d_j|^{1-\la}|d_i|^\la+|d_i|^{1-\la}|d_j|^\la}{|d_i|+|d_j|},
\max_{\theta_i\neq \theta_j }|\la\big(e^{i(\theta_j-\theta_i)}-1\big)+1| \right\}<1.
$$
  \item If $U\in\matu$ satisfies $N=UDU^*$, then $\sub{\cE}{N,\la}^s=U(\sub{\cE}{D,\la}^s) U^*$.
  \end{enumerate}
  In particular, the map $\orbu{D} \ni N \mapsto \sub{\cE}{N,\la}^s$ is smooth.  
  This fact can be formulated in terms of the projections $\sub{P}{N,\la} $ onto 
  $\sub{\cE}{N,\la}^s$ parallel to $\sub{T}{N}\orbu{D}$, $N \in \orbu{D}$. 
    \end{thm}
\proof See Section \ref{la prueba}.\QED

\begin{cor}\label{entorno y algo mas}
Let $D= \mbox{\rm diag}(d_1,\ldots,d_r)  \in \mat$ be an invertible diagonal matrix and $\la \in (0,1)$.
Let $\sub{\cE}{N, \la}^s$ and $k_D$ as in Theorem \ref{EL TEO}. Then, in 
$\orb{D}$ there exists a $\Delta_\la$-invariant $C^\infty$-pre-lamination 
$\{\ewe_{N, \la} \}_{N\in\orbu{D}}$ of class $C^\infty$ such that, for every $N \in \orbu{D}$, 
\begin{enumerate}
	\item $\ewe_{N, \la} $ is a $C^\infty$ submanifold of $\orb{D}$.
	\item $\sub{T}{N}\ewe_{N, \la} =\sub{\cE}{N, \la}^s\,$.
	\item   If \  $k_D <\rho<1$, then $\dist(\alulit{n}{T}-N)\leq \dist(T, N)  \rho^n$, for 
	every $T\in \ewe_{N, \la} \,$.
	 
	\item If $N_1\neq N_2$ then $\ewe_{N_1\, , \la}\cap\ewe_{N_2\, , \la}=\varnothing$.
  \item There exists an open subset $\ewe(D)$ of $\orb{D}$ such that 
  \ben
  \item [a. ]
  $\orbu{D} \inc  \ewe(D) \inc \displaystyle \bigcup_{N\in\orbu{D}}\ewe_{N, \la} $, and 
  \item [b. ] The map $p:\ewe(D)\to \orbu{D}$, defined by 
  $p(T)=N$ if $\,T\in\ewe_{N, \la} $, is of 
  class $C^\infty$.
  \een
\end{enumerate}
\end{cor}
\proof By Theorem \ref{EL TEO}, for every $k_D< \rho < 1$, 
$\orbu{D}$ is $\rho$-pseudo hyperbolic for $\Delta_\la$
(see Definition \ref{pseudohyperbolic}), and it consists of fixed points. Thus, by 
Theorem \ref{teorema 5.5}, we get a $C^\infty$ and $\Delta_\la$-invariant pre-lamination of class $C^\infty$, 
$\{\ewe_{N, \la} \}_{N\in\orbu{D}}$ which satisfies all the properties 
of our statement. 
\QED

\subsection{Convergence for fixed $\la$.}

\noi
Using the previous results, we can apply exactly the same techniques as in our previous work 
\cite{[APS]}, and to obtain for every $\la \in (0,1)$ the same results about $\alul{\cdot}$ as those 
obtained for the classical Aluthge transform $\alum{\cdot}$. We state these properties in the following
Theorem. The basic idea is to apply Proposition \ref{puntos limites normales} 
in order to assure that the iterations 
go into the open set $\ewe(D)$, where the smooth projection $p$ can be used. 
Although the proof of this theorem is omitted to avoid repetitions, Proposition \ref{SL} 
below gives a detailed proof.

\begin{thm} \label{tutti} \rm 
Let $\la \in (0,1)$:
\ben
\item Given a {\bf diagonalizable} matrix $T\in\mat$, the sequence $\{\alulit{n}{T}\}_{n\in\N}$ 
converges and its limit will be denoted by $\Delta_\la^\infty (T)$.
\item Let $D\in \mat$ be diagonal.  
Then the sequence $\{\Delta_\la ^n\}_{n\in\N}$, 
restricted to the similarity orbit $\orb{D}$,  converges uniformly on 
compact sets to  the map  $\Delta_\la^\infty:\orb{D}\to\orbu{D}$, which is of class $C^\infty$. 
In particular, $\Delta_\la^\infty$ is a $C^\infty$ retraction from $\orb{D}$ onto $\orbu{D}$.
\een
\end{thm}

\begin{rem}\label{largas} \rm 
Let $D\in \mat$ be diagonal.  For every $N \in \orbu{D}$ and $\la \in (0,1)$, denote 
$$
\ewe_{N, \la} ^+ = \big\{ \ T \in \orb{D} : \alulit{\infty}{T} = N \ \big\} \ .
$$
Since $T \mapsto \alulit{\infty}{T}$ is a $C^\infty$ retraction from $\orb{D}$ 
onto $\orbu{D}$, each $\ewe_{N, \la} ^+  $ is a $C^\infty$ submanifold of $\orb{D}$;
 and $\orb{D}$ stands as  the (disjoint) union of these sheets. On the other hand, 
 the submanifolds  $\ewe_{N, \la} ^+ $ are prolongations of the sheets 
$\ewe_{N, \la}  $ of Corollary \ref{entorno y algo mas}. Indeed, 
for every $N \in \orbu{D}$ and $\la \in (0,1)$, 
$$
\ewe_{N, \la}\inc \ewe_{N, \la} ^+ \peso{and} 
\ewe (D) \cap \ewe_{N, \la} = \ewe (D) \cap \ewe_{N, \la} ^+  \ ,
$$
by Corollary \ref{entorno y algo mas}. Then 
$\sub{T}{N}\ewe_{N, \la} ^+ =\sub{T}{N}\ewe_{N, \la}  =\sub{\cE}{N,\la}^s\,$.  \EOE
\end{rem}

\section{Regularity properties of $\alulit{\infty}{T} $}
As in section 3, we fix $D = \mbox{\rm diag}(d_1, \dots , d_r) \in \matinv$.
Observe that, using the continuity of $\alul{\cdot}$ with respect to $\la$ 
(Thm. \ref {continuidad}), and the 
fact that the convergence of its iterations to the limit map $\alulit{\infty}{\cdot}$ 
is uniform on compact subsets of $\orb{D}$ (for each $\la$), one can show that the map
$$
(0,1) \times \orb{D} \ni (\la , T) \ \longmapsto \ \alulit{\infty}{T} 
$$
is continuous. The purpose of this section is to prove that this map is smooth. 
Firstly we shall analyze the invertible case. Then, in subsection 4.2 we shall see 
that the invertibility hypothesis on $D$ can be dropped. 

\subsection{On the orbit $\orb{D}$.}
 Denote by 
$\orbl{D} = (0,1) \times \orb{D}$
and $\orbul{D} = (0,1) \times \orbu{D}$. Consider the map
\begin{equation}\label{laD}
\Delta : \orbl {D} \to \orbl{D} \ \ \mbox{given by} \ \ 
\Delta (\la , T) = (\la , \alul{T}\,  ) \ , \ \ (\la , T) \in \orbl{D} \ .
\end{equation}
\begin{rem}\label{previa} \rm
Using Theorems \ref{continuidad}, \ref{EL TEO} and \ref{tutti} and Proposition \ref{puntos limites normales}, 
one can deduce easily  the following properties:
\ben
\item $\orbl{D}$ is a $C^\infty$ manifold, and $\orbul{D}$ is a submanifold of $\orbl{D}$. 

\item The map $\Delta $ is of class $C^\infty$. 

\item The submanifold $\orbul{D}$ coincides with the set of  all fixed points of $\Delta$. 

\item For every $(\la , T) \in \orbl{D} $, it holds that 
\begin{equation}\label{lainf}
\Delta ^n  (\la , T) = (\la , \alulit{n}{T}\,  ) \conv (\la , \alulit{\infty}{T}\,  ) \ .
\end{equation}
\item For every $N \in \orbu{D}$ and $\la \in (0,1)$, one can describe the tangent spaces as
$$
\sub{T}{(\la , N)}\orbul{D}=\sub{T}{N}\orbu{D} \oplus\R  \ , \peso{and}
$$
$$
\sub{T}{(\la , N)}\orbl{D}= \sub{T}{N}\orb{D}\oplus\R  = \sub{\cE}{N,\la}^s\oplus 
\sub{T}{N}\orbu{D} \oplus \R 
=
\sub{\cE}{N,\la}^s\oplus  \sub{T}{(\la , N)}\orbul{D}  \ ,
$$
where the spaces $\sub{\cE}{N,\la}^s$ are those of Theorem \ref{EL TEO}. 
\item By Remark \ref{pnl6}, the projections $ \sub{Q}{N,\la} = \sub{P}{N,\la} + P_\R $ onto 
$\sub{\cE}{N,\la}^s$ parallel 
to $\sub{T}{(\la , N)}\orbul{D}\,$, given by the above decomposition, satisfy  that the map 
$$
(0, 1) \times \orbu{D} \ni (\la , N ) \longmapsto  
\sub{Q}{N,\la} 
$$
is of class $C^\infty$.
\EOE
\een
\end{rem}
\begin{prop} \label{SL}
The iterates $\Delta ^n$ of the map $\Delta : \orbl {D} \to \orbl{D}$ defined in Eq. \eqref{laD} satisfy 
$$
\Delta^n (\la , T) \conv \Delta^\infty (\la , T) \peso{for every} 
(\la , T) \in \orbl{D} \ , 
$$
where $\Delta^\infty : \orbl {D} \to \orbul{D}$ is a $C^\infty$ retraction. 
\end{prop}
\proof
Fix $(\la_0 , T_0 ) \in \orbl{D}$. Let $0<\la_1 < \la _0 <\la_2<1 $, 
and consider the submaniflods  
$$
\orbla{D} = (\la_1 \, , \, \la_2) \times \orb{D} \inc \orbl{D} 
\peso {and} \orbula{D} = (\la_1 \, , \, \la_2)\times \orbu{D} \ .
$$ 
Observe that $\orbla{D}$ is open in $\orbl{D}$,  it is   $\Delta$ invariant, and its fixed points 
coincide with $\orbula{D}$, which is also open in $\orbul{D}$. 
Fix $\rho \in (0,1)$ such that the constants (of Theorem \ref{EL TEO})  
$\sub{k}{D, \la} < \rho$ for every $\la \in (\la_1 \, , \, \la_2)$. 
Observe that, for every  $(\la , N) \in \orbula{D}$, 
$$ 
T_{(\la , N)}\,\Delta \Big|_{\sub{\cE}{N,\la}^s} = 
 \ T_N\,\Delta_\la\Big|_{\sub{\cE}{N,\la}^s} \inc \sub{\cE}{N,\la}^s \inc 
 \sub{T}{(\la , N)}\orbl{D} \ .
$$ 
Using this fact, and items 5 and 6  of  Remark \ref{previa}, one can assure that 
$\orbula{D}$ is $\rho$-pseudo hyperbolic for $\Delta : \orbla{D} \to \orbla{D} $ 
(see Definition \ref{pseudohyperbolic}) consisting of fixed points. Thus, by 
Theorem \ref{teorema 5.5}, we get a $C^\infty$ and $\Delta$-invariant pre-lamination 
$\{\ewe_{(\la , N)}\}_{(\la , N)\in\orbula{D}}$  of class $C^\infty$ such that, for every 
$(\la , N)\in\orbula{D}$, 
\ben
	\item $\ewe_{(\la , N)}$ is a $C^\infty$ submanifold of $\orbla{D}$.
	\item $\sub{T}{(\la , N)} \ewe_{(\la , N)} =\sub{\cE}{(\la , N)}^s\,$.
	\item    $\dist(\aluit{n}{\mu, T }-(\la , N)\,)\leq \rho^n \dist(\, (\mu, T), (\la, N)\, )  $, for 
	every $(\mu, T ) \in \ewe_{(\la , N)}\,$. 
Observe that this implies that,  
\beq\label{mismo la}
\peso{if $(\mu, T ) \in \ewe_{(\la , N)}\,$,  then $\mu = \la$ .}
\end{equation} 	 
	\item If $(\mu_1 , N_1 ) \neq (\mu_2 , N_2 ) $ then $\ewe_{(\mu_1 , N_1 ) }\cap
	\ewe_{(\mu_2 , N_2 )}=\varnothing$.
  \item There exists an open subset $\ewe(D)$ of $\orbla{D}$ such that 
  \ben
  \item [a. ]
  $\orbula{D} \inc  \ewe(D) \inc \displaystyle \bigcup_{(\la, N) \in\orbu{D}}\ewe_{(\la, N)}\ $. 
  
  \item [b. ] The map $p:\ewe(D)\to \orbula{D}$, defined by 
  $$
  p(\la, T)=(\la , N) \peso{ if } (\la , T) \in\ewe_{(\la, N)} \ ,
  $$ 
  is well defined and of   class $C^\infty$.
  \een  
\een
By item 2 of  Remark \ref{previa}, the map $\Delta$ and its iterations are $C^\infty$ functions. 
Hence, there exist an open set $\cU \inc \orbla{D}$ and $k \in \N$  such that 
$(\la_0 , T_0 ) \in \cU $ and $\Delta ^k (\cU )\inc \ewe(D)$. 
By the regularity the projection $p$, one can deduce that the map 
$
\orbla {D} \ni (\la , T) \ \longmapsto \ p (\Delta^k (\la , T) \,)  
$ 
is of  class $C^\infty$ on $\cU$. Now, 
if $\Delta^k (\la , T)  \in  \ewe_{(\la , N)}\,$ then 
$$
\displaystyle p (\Delta^k (\la , T) \,)   = N = 
\lim_{n\to \infty} \Delta ^n (\la , T) \ ,
$$
by item 3.  Hence $\Delta ^\infty$ is well defined and of   class $C^\infty$.    \QED

\begin{thm} \label{con la}
Let $D\in \diagr$ be invertible. Then 
the map 
$$
\Delta_\infty : (0,1) \times \orb{D} \to \orbu{D} \peso{given by}
\Delta_\infty (\la , T) = \alulit{\infty}{T}  \ , \ \ (\la , T) \in \orbl{D} \ 
$$
is of class $C^\infty$. 
\end{thm}
\proof  Let $\Delta^\infty : \orbl {D} \to \orbul{D}$ be the $C^\infty$ retraction of Proposition 
\ref{SL}. By Eq. \eqref{lainf}, 
$$
  \Delta^\infty (\la , T) = ( \la , \Delta _\infty (\la , T) \ ) \peso{for every } 
 (\la , T) \in \orbla {D} \ .
 $$
 This implies that also the map $\Delta _\infty$ is of class $C^\infty$.  \QED
 
 \begin{rem}\label{largas2} \rm 
Observe that Proposition \ref {SL} and Eq. \eqref{lainf}, before taking the limit, 
also show Theorem \ref{tutti}, whose proof was omitted. 
On the other hand, they also shows that the map 
$
\Delta^\infty : \orbl {D} \to \orbul{D}$ given by 
$\Delta^\infty (\la , T ) = (\la , \alulit{\infty}{T} ) $
is a $C^\infty$ retraction from $\orbl{D}$ onto $\orbul{D}$. Therefore, for every $(\la , N) \in 
\orbul{D}$, the set $\ewe_{(\la, N)}^+ = (\Delta^\infty )\inv (\la , N)$ is a smooth submanifold 
of $\orbl{D}$. Observe that, by Remark \ref{largas}, Eq. \eqref{lainf}
 and \eqref{mismo la}, we have that 
$\ewe_{(\la, N)}^+ = \{\la\} \times \ewe_{\la, N}^+$ and 
$$
\ewe_{(\la, N)} = \ewe(D) \cap \ewe_{(\la, N)}^+ = \ewe(D) \cap \Big( \{\la\} \times \ewe_{\la, N}^+ \Big) 
= \{\la\} \times \ewe_{\la, N} \ ,
$$
where $\ewe_{\la, N}^+ $ are the sheets described in Remark \ref{largas}, and 
$\ewe_{\la, N}$ are those sheets which appear in Corollary \ref{entorno y algo mas}. 
This can be described as the fact that ``the submanifolds $\ewe_{\la, N}$ move smoothly
with $\la$". 
\EOE
\end{rem}

\subsection{The not invertible case}

As in section 3, all the results of this section still hold if the diagonal 
matrix $D\in \mat$ is not invertible.  Indeed, 
suppose that $\rk{D} = s<r$, and fix $T\in\orb{D}$ and $\la \in (0,1)$. 
By Eq. \eqref{tuno}, 
$$
\alul{T}=\begin{pmatrix}
T_1&0\\
0  &0
\end{pmatrix}\begin{array}{l}
\ker T ^\bot\\
\ker T 
\end{array}\, ,
$$
where $T_1$ is invertible and diagonalizable on $\ker T ^\bot$ . The same happens 
for every matrix $S\in\orb{D}$. 
Denote by $ P = P_{\ker T} $ and $\orbu{P} = \{ UPU^* : U \in \matu\} $ its unitary orbit. 
Consider the map 
$$ Q : \orb{D} \to \orbu{P} \peso{given by} Q(S) = P_{\ker S} = P_{\ker \alul{S}} \ , \quad 
S\in\orb{D} \ .
$$
This map takes values in $\orbu{P}$ because $\rk{Q(S)} = \rk{P}=s$ for 
every $S \in \orb{D}$. Also, $Q$ is of class $C^\infty$, 
since $S \mapsto \alul{S}$ is smooth, and 
$\alul{S} \mapsto P_{\ker \alul{S}}$ is polynomial. 
By Proposition \ref{son variedades}, there exist an open set 
$\cW \inc \orbu{P}$ which contains $P$, 
and a $C^\infty$ local cross section $\sigma : \cW \to \matu$, in the sense that 
$\sigma (R) P \sigma (R)^* = R$ for every $R \in \cW$. 
Let $\cV \inc \orb{D}$ be an open set such that $T\in \cV$ and $Q(\cV ) \inc \cW$. Denote by 
$\eta = \sigma \circ Q  : \cV \to \matu$. Then $\eta$ is also
of class $C^\infty$.  So, for every $\la \in (0,1)$ and $S \in \cV$, 
there exists $\gamma(S, \la)\in L(\ker T\orto )$ such that 
$$
\eta(S)^*  \, \alul{S}\, \eta(S)= 
\begin{pmatrix}
\gamma(S, \la)&0\\
0  &0
\end{pmatrix}\begin{array}{l}
\ker T ^\bot\\
\ker T  
\end{array} \ \ \mbox{and} \ \ \gamma(S, \la)\in\orb{T_1} \inc L(\ker T \orto ) \ .
$$
Therefore, using that $\alul{UTU^*}=U\alul{T}U^*$ for 
every  $U\in \matu$,  we obtain
$$
 \alulit{\infty}{S}= \sub{Ad}{\eta(S)}\circ \begin{pmatrix}
\alulit{\infty}{\gamma(S, \la)\,}&0\\
0  &0
\end{pmatrix}\begin{array}{l}
\ker T ^\bot\\
\ker T 
\end{array}\  ,\peso{for every  $S\in\ewe$ .}
$$
Then the regularity of $(\la , S) \mapsto \alulit{\infty}{S}$ can 
be deduced from the regularity of the maps 
$(0,1) \times\cV\ni (\la ,S) \mapsto \gamma(S, \la)$ and 
$(0,1) \times \orb{T_1} \ni (\la , A) \mapsto 
\alulit{\infty}{A}$. 
Hence, the reduction to the invertible case is proved.
 
\subsection{Different eigenvalues.}

Let  $\diagdis$ be the set of diagonalizable and invertible matrices in $\mat$ 
with $r$ different eigenvalues (i.e. every eigenvalue has algebraic multiplicity equal 
to one). Observe that $\diagdis$ is an open dense subset of $\mat$ and it is invariant 
by the  Aluthge transform.

\begin {thm}\label{en diagdis}
The map $\Delta_\infty : (0,1) \times \diagdis \to \orbu{D}$ given by
$$
\Delta_\infty (\la , T) = \alulit{\infty}{T}  \ , \ \ (\la , T) \in (0,1) \times \diagdis \ 
$$
is of class $C^\infty$. 
\end{thm}
\proof 
It follows from  a straightforward combination of the techniques of section 3.2 of \cite {[APS]} and 
those of the previous sections. We omit the details. \QED

\subsection{The map $\la \mapsto \alulit{\infty}{T}$ for fixed $T$.}

\begin{defn}\label{RT} \rm 
Let  $T \in \mat$ be a diagonalizable. We call 
$R_T : (0,1) \to \mat$ the map given by $R_T (\la ) = \alulit{\infty}{T}$, for $\la \in (0,1)$. 
\end{defn}

\noi
The following question arises naturally: Is the map 
$R_T$ constant for every  diagonalizable matrix $T \in \mat$?. 
Numerical examples  shows that the question has a negative answer in general. 
Indeed, taking the $2\times 2$ matrix $T= \bm{rl} 3&0\\-2&1\em$, numerical computations show that 
\begin{align*}
R_T(0.3) \cong \bm{cc} 2.2273 &0.97380\\0.97380& 1.7726 \em 
 \ \mbox{and} \ \ 
 R_T(0.7)  \cong \bm{cc} 1.37162 &-0.77790\\-0.77790& 2.62838 \em \ .
\end{align*}
Nevertheless, it was proved in \cite[Thm. 4.9]{[AMS]} that, if 
$D=\mbox{\rm diag}(d_1,d_2)\in \mathcal{M}_2 (\C)$ and $|d_1|=|d_2|$, then $R_T$ 
is constant for every $T\in\orb{D}$. Our next result shows that, for a 
diagonal matrix $D=\mbox{\rm diag}(d_1,\ldots,d_r)\in\mat$, this may 
happen only if $|d_1|=\ldots=|d_r|$.

\begin{prop} \label{nocon} Let $D 
\in \mat$ be diagonal and invertible. If $D$ has two eigenvalues with distinct moduli, 
then there exists $T \in \orb{D}$ such that the map $R_T$  is not constant. 
\end{prop}

\noi The proof of Proposition \ref{nocon} follows directly from the next two Lemmas. But let us 
make first some comments about this problem:

\begin{rem}\label{uno no}\rm
Despite Proposition \ref{nocon}, given $D = \mbox{\rm diag} \, ({d_1 \, , \dots , d_r}) \in \mat$ 
with $|d_i|\neq |d_j|$ 
for some $i,j\in\{1,\ldots,r\}$, the similarity orbit $\orb{D}$ 
may contain (non normal) matrices $T$ for which the map $R_T$ 
is constant. In fact, consider the following example: 
let $D = \mbox{\rm diag} \, ({1 , -1 , 2}) \in \cM_3(\C )$, 
take  $D_1 = \mbox{\rm diag} \, ({1 , -1 })$, $T_1 \in \orb{D_1}$ and 
 $T = \bm{cc} T_1 & 0\\ 0&2 \em\in \orb{D}$, then 
 $$
 R_T(\la ) =  \alulit{\infty}{T} = \bm{cc} \alulit{\infty}{T_1}& 0\\ 0&2 \em 
 = \bm{cc} R_{T_1} (\la )& 0\\ 0&2 \em \ , 
 \ \mbox{for every} \ \la \in (0,1)\ .
 $$
 But the map $R_{T_1}$ is constant by  \cite[Thm 4.9]{[AMS]} 
 (see also Proposition \ref{con dos} below). Observe that this example 
 includes several cases where $T \notin \orbu{D}$ (otherwise, the map would be trivially constant). 
  On the other hand, Proposition \ref{nocon} is not longer true if $D $ 
 is not invertible. This fact can be immediately tested by taking $D$ with 
 $\spec{D} = \{0,1\}$ (i.e., if $D$ is a projection). 
 \EOE
\end{rem}

\begin{lem} \label{constorb} Let $D \in \mat$ be diagonal and invertible. If 
the map $R_T$ is constant for every $T \in \orb{D}$,  
then the distribution of subspaces 
$$
(0,1) \ni \la \longmapsto \sub{\cE}{N,\la}^s 
$$
given by Theorem \ref{EL TEO} must be constant for every $N \in \orbu{D}$. 
\end{lem}
\proof Using the notations of 
Remark \ref{largas} and Corollary \ref{entorno y algo mas}, 
each submanifold $\ewe_{N, \la}^+ $ consist of those matrices $T$ 
such that $R_T(\la) = \alulit{\infty}{T} = N$. But if all the maps $\alulit{\infty} {\cdot}$ 
are the same on $\orb{D}$, 
then the submanifolds $\ewe_{N, \la}^+ $ must agree for different lambdas. 
Finally observe that, by  Remark \ref{largas} and 
Corollary \ref{entorno y algo mas}, $\sub{\cE}{N,\la}^s $ is the tangent space of $\ewe_{N, \la}^+ $ 
at $N$, for every $N \in \orbu{D}$ and every $\la \in (0,1)$. \QED

\begin{lem} \label{no varian} Let $D = \mbox{\rm diag} \, ({d_1 \, , \dots , d_n}) \in \mat$ 
be diagonal and invertible. Then the following conditions are equivalent:
\ben
\item The distribution of subspaces $(0,1) \ni \la \longmapsto \sub{\cE}{N,\la}^s $ 
given by Theorem \ref{EL TEO} is constant for every $N \in \orbu{D}$. 
\item  $|d_i | = |d_j| $ for every $1\le i,\, j \le r$. 
\een
Moreover, in this case, $\sub{\cE}{N,\la}^s = \big[ \sub{T}{D}\orbu{D} \big]^\bot$ for every $\la \in (0,1)$. 
\end{lem}
\proof Since the proof uses several results and notation from section 5, we postpone it 
until that section. See Remark \ref{varia}. 
\QED

\medskip

\noi Now, the natural question is, what happens if $D$ has all its 
eigenvalues of the same modulus?. We first study a particular case: 

\begin{prop} \label{con dos}
 Let $D\in \mat$ be diagonal such that $\sigma (D) = \{ d_1 \, , d_2\}$ with  $|d_1 | = |d_2 |$ . 
 Then the map $R_{T}$ is constant for every $T\in \orb{D}$.
 \end{prop}
 \proof  Let $T \in \orb{D}$. Denote  $\ese _i = \ker (T-d_i I)$ for $i = 1,2$, 
  and by $Q$ the oblique projection
 onto $\ese_1$ given by the decomposition $\cene = \ese_1 \oplus \ese_2\,$. Then $T = d_1 Q + 
 d_2 (I-Q)$. 
Assume that  $k = \dim \ese_2  \ge n = \dim \ese_1  \,$.  
In matrix terms, there exists $A \in L(\ese_1\orto \, , \, \ese_1 )$ such that  
 $$
Q = \bm {cc} I & A \\ 0 & 0 \em \barr{c} \ese_1 \\\ese_1\orto \earr   \peso {and} 
T = \bm {cc} d_1 \, I & (d_1  - d_2) A  \\ 0 & d_2 I  \em  \barr{c} \ese_1 \\\ese_1\orto \earr \ .
$$
We can assume that 
$\ese_1 = \{ x \in \cene : x_i = 0 $ for $i > n \}$, by a unitary conjugation, 
which commutes with $\Delta_\la^\infty$. In this case, $A \in \cM_{n,k}(\C) $. 
By the decomposition in singular values of $B= (d_1  - d_2) A $, 
there exist $U \in \cU(n) $ and $V \in \cU(k)$ such 
that $UBV^* = \Sigma (B)$, where $\Sigma (B) = (\Sigma_0 (B) , 0) \in \cM_{n,k}(\R) $ 
with $\Sigma_0 (B) = \diag{s_1(B) , \dots , s_n(B)\,} \in \cM_n (\R)$,  
and  we add a $n\times (k-n)$ block of zeros on the right. 
If $W = U \oplus V \in \matu$, then 
$$
WTW^* = \bm {cc} d_1 \, I & \Sigma (B) \\ 0 & d_2 I  \em  \barr{c} \C^n  \\ \C^k \earr 
= \bm {ccc} d_1 \, I_n & \Sigma_0 (B)& 0  \\ 0 & d_2 \, I_n &0 \\
0 & 0 & d_2 \, I_{k-n}  \em  \barr{l} \C^n  \\\C^n  \\ \C^{k-n} \earr 
\ . 
$$ 
Moreover, there exists  a permutation matrix $S \in \matu\,$ 
which rearranges the entries of $WTW^*$ in such a way that
$$
M = SWTW^*S^* =  \bigoplus _{i = 1} ^n \bm {cc} d_1& s_i (B) \\ 0 & d_2   \em 
\oplus d_2 \, I_{k-n} \ .
$$
As before, it suffices to show that the map $R_M$ is constant. But now 
Proposition \ref{algeo} assures that 
$$
R_M(\la ) = \alulit{\infty}{M} = \bigoplus _{i = 1} ^n \Delta_\la^\infty 
\bm {cc} d_1& s_i (B) \\ 0 & d_2   \em 
\oplus d_2 \, I_{k-n} \peso{for every $\la \in (0,1)$ .} 
$$
Finally, it was proved in \cite[Thm. 4.9]{[AMS]} that each map 
$\la \mapsto  \Delta_\la^\infty \bm {cc} d_1& s_i (B) \\ 0 & d_2   \em $ is constant. 
\QED

\begin{rem}\label{largo largo largo} \rm
The case $D = D^*$  in Proposition \ref{con dos}, is particularly 
interesting, because in this case all iterated $\la$-Alutge transforms (and the limit) can 
be explicitly computed for every $E \in \orb{D}$, even in the infinite dimensional case:
Let $\hil$ be a Hilbert space and	$D\in \opsa$ such that $\sigma(D) = \{ 1, -1\}$. 
Observe that if $E\in \orb{D}$, then $E^2 = E$. 
The geometry of $\orb{D}$ has been widely studied (see \cite{[PR2]}, \cite{[CPR2]} and \cite{[CPR3]}). 
Given $E \in \orb{D}$, observe that, if $L = |E | $, then $|E^* | = L\inv$. 
Therefore, if $E = R L$ is the polar decomposition of $E$, with $R \in \matu$, then
$RL = E = E\inv = L\inv R^*$. But the right  polar decomposition is  $E = |E^*| R$ (with the same $R$), 
so that $R = R\inv = R^*$ is a unitary reflection and $RL = L\inv R$. Moreover, 
$R \in \orbu{D}$, because $E = RL = L\mrai R L\rai$.

The map $E \mapsto R$ is the retraction $p : \orb{D} \to \orbu{D}$ deeply studied  in \cite {[PR3]}. 
We claim that $p = \alus{\cdot} = \alulit{\infty}{\cdot} $ for every $\la \in (0,1)$. 
Indeed, observe that 
$$
\alus{E} = L\rai R L\rai = L\rai L\mrai R = R 
$$ 
and, for every $\la \in (0,1)$,  $
\alul{E} = L^\la R L^{1-\la} = RL^{1-2\la}$, so that $ 
\alulit{2}{E} = RL^{(1-2\la)^2}  \ ,
$
since $\alul{E}= RL^{1-2\la} $ is the polar decomposition of $\alul{E}$. Inductively,  
for every $n \in \N$,  
$$
\alulit{n}{E} = RL^{(1-2\la)^n} \conv R = \alulit{\infty}{E} \ ,
$$
since  $|1-2\la|<1$. Then the map $\la \mapsto R_E(\la ) = \alulit{\infty}{E}$ is constant, but the rate 
of convergence is very different for each $\la$, being slower when $\la$ tends to $0$ or $1$. 
On the other hand, the sheets $\ewe_{R, \la}^+ $ of Remark \ref{largas} can be 
characterized as 
$$
\ewe_{R, \la}^+ = \{\  RM \ : \ M \in \glh^+   \peso{and}  RM = M\inv R\ \} \ .
$$ 
The geometry of these hyperbolic manifolds is also deeply studied in \cite{[CPR3]}. \EOE
\end{rem}

\smallskip
\noi
If $D$ has all its eigenvalues of the same modulus and $\spec{D}$ has more than two elements, 
we do not have an answer to the above question, but we have made several computational experiments.
In all the tested examples,  
the map $R_{T}$  fails 
to be constant for some $T$ in the orbit, even if $D$ satisfies some algebraic condition 
such as $D^{3}=I$. This suggests the following conjecture\,: 

\begin{conj}\label{con} \rm
 Let $D\in \mat$ be diagonal and invertible. Then the map $R_{T}$ is constant for every $T\in \orb{D}$ \iiff   $\sigma (D) = \{ d_1 \, , d_2\}$ with  $|d_1 | = |d_2 |$ . \EOE
\end{conj}

\begin{rem} \rm
Observe that it would be  sufficient to consider the $3 \times 3$
case, because we can use a similar reduction to the one used in Remark \ref{uno no}. 
\EOE
\end{rem}

\begin{exa}\rm
Let $a,\,b,\, c\in\R^+$ such that $abc=1$, and let $T$ be the $3\times 3$ matrix defined by
$$
T=\begin{pmatrix}
0&0&1\\
1&0&0\\
0&1&0
\end{pmatrix}\cdot
\begin{pmatrix}
a&0&0\\
0&b&0\\
0&0&c
\end{pmatrix}=U|T| \ .
$$
Computing its spectrum, one shows that $T\in \orb{U}$. 
On the other hand, since $U$ is a permutation matrix, for every diagonal matrix $D\in \cM_3(\C)$ both  
$UDU^*$ and $U^*DU$ are also diagonal matrices. In particular, $|T^*|=U|T|U^*$ is diagonal and  
commutes with $|T|$. Let $\la\in(0,1)$. Then
\begin{align*}
\alul{T}=\alulf{T} =U\big(U^*(|T|^\la|T^*|^{1-\la})U\big)=U\big|\alul{T}\big| \ ,
\end{align*}
where the last equality holds by the uniqueness of the polar decomposition. 
Note that $\big|\alul{T}\big|$ is diagonal with 
$\det \big|\alul{T}\big| =1$. An inductive argument shows that
$\alulit{n}{T}=UD_n $ for every $ n\in\N$, 
where each $D_n$ is a positive diagonal matrix.
As $T\in \orb{U}$,  then $\alulit{\infty}{T} \in \orbu{U} \inc \matu$. 
Hence $D_n \conv I$ and  $\alulit{\infty}{T}=U$. The same happens for any 
$\la\in(0,1)$, therefore, the function $R_T(\la)$ is constant.
This example does not contradicts the Conjecture, because 
all the matrices studied satisfy that $|T|$ and $|T^*|$ commute, so that 
they are not dense in  $\orb{U}$. \EOE
\end{exa}

\section{The proof of Theorem \ref{EL TEO}.}\label{la prueba}
As in Section 3, in this section we fix an invertible diagonal matrix $D\in\mat$ whose diagonal 
entries are denoted by $(d_1,\ldots,d_n)$. For every $j\in\{1,\ldots,n\}$, let  
$d_j=e^{\,i\theta_j}|d_j|$ be the polar decomposition of $d_j$, where $\theta_j\in [0,2\pi]$.

\subsection{Matricial characterization of $T_N \Delta_\la $}

\begin{defn}\label{producto de Hadamard}\rm
Given  $A,B\in\mat$, $A\circ B$ denotes their Hadamard product, 
that is, if $A=(A_{ij})$ and $B=(B_{ij})$, then $(A\circ B)_{ij}=A_{ij} B_{ij}\,$. 
With respect to this product, each matrix $A\in\mat$ induces an operator $\had{A}$ 
on $\mat$ defined by $\had{A}(B)=A\circ B$, $B \in \mat$. 
\end{defn}

\begin{rem}\label{remark sobre hadamard}\rm
In what follows, we shall state several definitions and results taken from Section 
4 of our previous work \cite{[APS]}: 
\ben
\item By Eq. \eqref{tan1}, the tangent space $\sub{T}{D}\orb{D}$ consists on those matrices 
$X \in \mat$ such that $X_{ij} = 0$ if $d_i = d_j\,$. 
Then $\sub{T}{D}\orb{D}$  reduces the operator $\had{A}\,$, for every $A\in\mat$. 
This is the reason why, from now on, we shall consider all these operators as acting on 
$\sub{T}{D}\orb{D}$. Restricted in this way, it holds that
\[
\|\had{A}\|=\sup\{\|A\circ B\|_{_2}: \ B\in\sub{T}{D}\orb{D} \ \mbox{and}\  
\|B\|_{_2}=1\}=\max_{d_i\neq d_j} |A_{ij}| \ ,
\]
since $\had{A}$ is a $diagonal$ operator on the Hilbert space $(\sub{T}{D}\orb{D} \ , \ \| \cdot \|_{_2})$. 

\item Let $\preal$ and $\pim$ be the projections defined on $\sub{T}{D}\orb{D}$ by
$$
\preal(B)=\frac{B+B^*}{2} \peso{and} \pim(B)=\frac{B-B^*}{2} \quad, \quad B \in \sub{T}{D}\orb{D}\ .
$$
That is, $\preal$ (resp. $\pim$) is the restriction to $\sub{T}{D}\orb{D}$ of the  orthogonal projection onto the subspace of hermitian (resp. anti-hermitian) matrices.  
\item 
Observe that, for every $K \in \matah$ (i.e., such that $K^* = -K$) and $B\in\mat$ it holds that 
\begin {equation}\label{commu}
K\circ \preal(B)=\pim(K\circ B)  \peso{and} K\circ \pim(B)=\preal(K\circ B)\ .
\end{equation}
\item Denote by $\qd\,$ 
the orthogonal projection from $\sub{T}{D}\orb{D}$ onto $(\sub{T}{D}\orbu{D})^\bot$. 
 \item 
 \label{pruyucciun}
Let $J,K\in\mat$ be the matrices defined by
\begin{align*}
K_{ij}=
\begin{cases}
|d_j-d_i|\sgn(j-i)&\mbox{if $d_i\neq d_j$}\\
0                 &\mbox{if $d_i = d_j$}
\end{cases}
\ \ \mbox{and} 
\end{align*}
\begin{align*}
J_{ij}=\begin{cases}
(d_j-d_i)K_{ij}^{-1}&\mbox{if $d_i\neq d_j$}\\
1                   &\mbox{if $d_i = d_j$}
\end{cases},
\end{align*}
for $1\leq i,j\leq r$. Then
\begin{enumerate}
\item For every $A\in\mat$, $AD-DA=J\circ K\circ A$.
\item  It holds that $\qd \ = \had{J}\pim\had{J}^{-1}$.
	
\item If $H \in \matsa$ (i.e., if $H^*=H$), then 	$\qd \had{H}=\had{H}\qd\,$.
\end{enumerate}

\item 
\label{la matriz}\rm
Let $\la \in (0,1)$ and $N \in \orbu{D}$ and let $\qn$ be the orthogonal projection from 
$\sub{T}{N}\orb{D}$ onto $\big(\sub{T}{N}\orbu{D}\big)^\bot$. 
Then $\dersin{\Delta_\la}{N}$ has the following $2\times 2$ 
matrix decomposition
\begin{equation}
\label{derivada en dos por dos}
\dersin{\Delta_\la}{N}=
\begin{pmatrix}
\sub{A}{1N}(\la ) & 0\\
\sub{A}{2N}(\la ) & I   
\end{pmatrix}  \barr{r} \qn\\ I- \qn \earr \ ,
\end{equation}
because $\dersin{\Delta_\la}{N}$ acts as the identity on $\sub{T}{N}\orbu{D}$. 

\item 
\label{derivada de gamma en cero}
Let $A\in \mat$ and let $\gamma:\R\to\posop$ be the curve defined by 
$$
\gamma(t)=\big(e^{tA}De^{-tA}\big)^*\big(e^{tA}De^{-tA}\big)=e^{-tA^*}D^*e^{tA^*}e^{tA}De^{-tA}.
$$
If $R$, $T^+$ and $T^- \in \mat$ are defined by 
$R_{ij}=2\bar{d}_i d_j$, $T^+_{ij}=|d_i|^2+|d_j|^2$ and 
$T^-_{ij}=|d_j|^2-|d_i|^2$, $i,j \in \mathbb{I}_{r}\,$, then 
\begin{align}
\gamma'(0)=(R-T^+)\circ\preal(A)+T^-\circ \pim(A).\label{gamma prima}
\end{align}
\een
\end{rem}

\newcommand{\matreal}{\mathcal{M}_n(\mathbb{R})}

\medskip
\noi
The following classical result, proved by Daleki\u\i i \ and Kre\u\i n in 1951 
(\cite{[DK1]} and \cite{[DK2]})   will be useful in the sequel 
 (see also the book \cite{[HJ2]}). 

\begin{thm} \label{DK} \rm 
Let $I, J \inc \R$ be open intervals and let  $\gamma : I \to \matsa$  be a $C^1$ curve
such that $\sigma (\gamma (t)\,) \inc J$ for every $t \in I$. Let  $f : J \to \R$ be 
a $C^1$ map. Suppose that $\gamma (t_0) = \diag{a_1 , \dots , a_r}$
 for some  $t_0 \in I$. Then 
$$
(f \circ \gamma )' (t_0) = M_f \circ \gamma ' (t_0) \ , 
$$
where $M_f \in \matreal$ is defined by 
$$
(M_f)_{ij} =\begin{cases}
\frac{\displaystyle f(a_j)-f(a_i)}
{\displaystyle a_j-a_i} & \mbox{if $a_i \neq a_j$} \\ & \\
f'(a_i) & \mbox{if $a_i = a_j$} \end{cases}  \quad  , 
\peso{for} i\ , \ j \ \in \mathbb{I}_r \ .  \label{matriz emeefe}
$$
\end{thm}

\begin{cor}\label{derivada de la curva}
Let $A$ and $\gamma$ be as in \ref{derivada de gamma en cero}  of Remark \ref{remark sobre hadamard}, 
and let $\la\in(0,1)$. 
Then  $(\gamma^{\la /2})'(0)=M_{\la /2}\circ \gamma'(0)$, where 
$M_{\la /2}\in \mat $ is the matrix given by 
\begin{align}
(M_{\la /2})_{ij}=\begin{cases}
\frac{\displaystyle |d_j|^{\la}-|d_i|^{\la}}{\displaystyle |d_j|^2-|d_i|^2} 
& \mbox{if $|d_i| \neq |d_j|$} \\ & \\
\frac {\displaystyle \la}{2} \ |d_i|^{\la-2}& \mbox{if $|d_i| = |d_j|$}
\end{cases}  \quad , \quad i\ , \ j \ \in \mathbb{I}_r  \ . \label{matriz eme} 
\end{align}
\end{cor}
\proof Apply Theorem \ref{DK} to $f(t) = t^{^{\frac\la 2}}$. Use that 
$\gamma (0) = \diag {|d_1|^2 , \dots , |d_r|^2 }$.
\QED

\bigskip
\noi 
Recall that  $\qd$ denotes the orthogonal projection from $\sub{T}{D}\orb{D}$ 
onto $\big(\sub{T}{D}\orbu{D}\big)^\bot$.

\begin{prop}\label{uno}
Fix $\la \in (0,1)$. Then there exists a matrix $H(\la ) \in \mat$ such that 
$\qd \big(\dersin{\Delta_\la}{D}\big) \qd =\qd\ \had{\preal H(\la )} \ \qd$
and  
$$
(I-\qd) \big(\dersin{\Delta_\la}{D}\big) \qd = (I-\qd)\  \had{\pim H(\la )} \ \qd  \, .
$$
The entries of $H_1(\la ) = \preal H(\la )$ are the following: 
for every $1 \le i, j \le r$, denote by $A_{ij} = |d_j|^{2-\la}|d_i|^\la-|d_i|^{2-\la}|d_j|^\la$
and $B_{ij} = |d_j|^{1+\la}|d_i|^{1-\la} - |d_i|^{1+\la}|d_j|^{1-\la}$. Then  
\begin{align}
H_1(\la )_{ij}  
& = \begin{cases}
\frac{\displaystyle A_{ij} +\ e^{i(\theta_j-\theta_i)} B_{ij}}
{\displaystyle |d_j|^2-|d_i|^2} & \mbox{\small{ if $|d_i|\neq |d_j|$}}\\&\\
\la\big(e^{i(\theta_j-\theta_i)}-1\big)+1&\mbox{\small{if $|d_i|=|d_j|$}}
\end{cases} \quad .
\label{h1} 
\end{align}
\end{prop}

\noi The proof of this proposition follows the same steps as the proof of 
 of \cite[Prop. 4.1.5]{[APS]}, but using now Corollary \ref{derivada de la curva}, 
 and items \ref{pruyucciun} and \ref{derivada de gamma en cero}  of 
 Remark \ref{remark sobre hadamard}. For this reason we shall give only and 
 sketched version, pointed out the main differences and the technical 
 difficulties that appear when we loose the symmetry of the case $\la=1/2$.

\bdem[Sketch of proof]
Fix  $X=AD-DA \in T_D \orb{D}$, for some $A \in \mat$. Then
$$
\der{\Delta_\la}{D}{X}=\left.\frac{d}{dt}\alul{e^{tA}De^{-tA}}\right|_{t=0}.
$$
Let $\gamma(t)=\big(e^{tA}De^{-tA}\big)^*\big(e^{tA}De^{-tA}\big)=e^{-tA^*}D^*e^{tA^*}e^{tA}De^{-tA}$. 
In terms of $\gamma$, we can write the curve $\alul{e^{tA}De^{-tA}}$ in the following way
$$
\alul{e^{tA}De^{-tA}}=\gamma^{\la/2}(t)(e^{tA}De^{-tA})\gamma^{-\la/2}(t).
$$
Since $\gamma^{\la/2}\gamma^{-\la/2}=I$, then 
$(\gamma^{-\la/2})'(0)=-\gamma^{-\la/2}(0)\ (\gamma^{\la/2})'(0)\ \gamma^{-\la/2}(0)$. 
Using this identity, easy computations  show that 
\begin{align*}
\der{\Delta_\la}{D}{X}&
= \Big((\gamma^{\la/2})'(0)\ D - D\ (\gamma^{\la/2})'(0)\Big) |D|^{-\la} +|D|^{\la}(AD-DA) |D|^{-\la}.
\end{align*}

\noi If we define the matrices $L,N\in\mat$ by
$N_{ij}=|d_j|^{-\la}$ and $ L_{ij}=|d_i|^{\la}|d_j|^{-\la}$ 
and take $J,K\in\mat$ as in \ref{pruyucciun}  of Remark \ref{remark sobre hadamard}. Then
\begin{align*}
\der{\Delta_\la}{D}{X}
&=N\circ( J\circ K \circ (\gamma^{\la/2})'(0))+L\circ(J\circ K\circ A).
\end{align*}
Using Corollary \ref{derivada de la curva} and  \ref{derivada de gamma en cero}  of 
Remark \ref{remark sobre hadamard},  we get
\begin{align*}
\der{\Delta_\la}{D}{AD-DA}&=N\circ J\circ K\circ M_{\la/2}\circ\Big[(R-T^+)\circ        
                       \preal(A)+T^-\circ\pim(A)\Big]\\
                       &+L\circ J\circ K\circ A.
\end{align*}
where $R$, $T^+$ and $T^-$ are the matrices defined in \ref{derivada de gamma en cero}  
of Remark \ref{remark sobre hadamard}, and 
$M_{\la/2}$ is the matrix defined in Eq. \eqref{matriz eme} of Corollary \ref{derivada de la curva}.  

\medskip 
\noi 
Now, we shall express $\der{\Delta_\la}{D}{AD-DA}$ in terms of $AD-DA=J\circ K\circ A$. 
Following the same steps as in Proposition 4.1.5 of \cite {[APS]}, we arrive to the formula 
$$
\der{\Delta_\la}{D}{AD-DA}=\Big(\, \had{H(\la )}\  \qd +(I-\qd)\,\Big)(AD-DA) \ ,
$$ 
where $H(\la )=M_{\la/2}\circ N\circ (R-T^+)+L$. 
Given $X\in R(\qd)$, then $\had{J}^{-1}(X)\in \matah $ and 
\begin{align*}
\qd\big(\dersin{\Delta_\la}{D}\big) \qd (X) 
&=\qd(H(\la )\circ X)=(\had{J}\pim\had{J}^{-1})(H(\la )\circ X) \\
&=J\circ\Big(\pim (H(\la )\circ \had{J}^{-1}X)\Big)\\
&=\frac12 \ J\circ\Big( H(\la )\circ \had{J}^{-1}(X)\,+\,H(\la )^*\circ \had{J}^{-1}(X)\Big) \\
&=J\circ \preal H(\la )\circ \had{J}^{-1}(X)= \preal H(\la )\circ X = \qd \had{\preal H(\la )} ( X ) \ .
\nonumber 
\end{align*}
Analogously, one shows that $(I-\qd)\big(\dersin{\Delta}{D}\big) \qd (X) = 
(I- \qd )  \had{\pim H(\la )} ( X )$. 
In order to prove Eq. \eqref{h1},  recall that $H(\la )=M_{\la/2}\circ N\circ (R-T^+)+L$. 
Hence, $H(\la)_{ij}=(M_{\la/2})_{ij}\Big(|d_j|^{-\la}\big(2\bar{d_i}d_j-(|d_i|^2+|d_j|^2)\big)\Big)+
|d_i|^{\la}|d_j|^{-\la}\,$. 
Suppose that $|d_i|\neq |d_j|$. Straightforward computations, using 
Corollary \ref{derivada de la curva},  show that 
\begin{align*}
\preal H (\la )_{ij}
&=\frac{1}{2}\left(\frac{2|d_j|^{2-\la}|d_i|^\la+2\bar{d_i}d_j-
2\bar{d_i}d_j|d_i|^\la|d_j|^{-\la}-|d_i|^2-|d_j|^2}{|d_j|^2-|d_i|^2}\right)+ \\
&\ \ \ \, \,  \frac{1}{2}\left(\frac{2|d_i|^{2-\la}|d_j|^\la+2\bar{d_i}d_j-
2\bar{d_i}d_j|d_j|^\la|d_i|^{-\la}-|d_j|^2-|d_i|^2}{|d_i|^2-|d_j|^2}\right)\\
&= \frac{A_{ij} \ +\ e^{i(\theta_j-\theta_i)}B_{ij}}{|d_j|^2-|d_i|^2} \ , 
\end{align*}
where $A_{ij}$ and $B_{ij}$ are those of the statement. 
If $|d_i|=|d_j|$, then 
$$
H(\la )_{ij}= \ \frac{\la}{2}|d_i|^{\la-2} \ \Big(|d_i|^{2-\la}2\big(e^{i(\theta_j-\theta_i)}-1\big)\Big)+1
=\la\big(e^{i(\theta_j-\theta_i)}-1\big)+1 \ .
$$  
So that $H(\la )_{ij} = \overline{H(\la )_{ji}} = 
\preal H(\la )_{ij}  = H_1(\la)_{ij} =  \la\big(e^{i(\theta_j-\theta_i)}-1\big)+1
$.
\edem

\begin{rem} \label{con a} \rm
Using the notations of Proposition \ref{uno}, let $H_2 (\la )= \pim H(\la )$. 
\ben
\item If $|d_i|=|d_j|$, as we observed at the end of the proof of Proposition \ref{uno}, $H(\la)_{ij}=\overline{H(\la)_{ji}}$. Hence,
$H_2(\la )_{ij} = 0$.
\item 
If $|d_i| \neq |d_j|$, denote $a = \frac{\displaystyle|d_j|}{\displaystyle|d_i|} \neq 1$ 
and $\alpha = e^{i(\theta_j-\theta_i)} $. 
Then one can show that 
$$
H_2(\la )_{ij} = (a-a\inv)\inv 
\ \Big[ \ a ^{1-\la } + a ^{\la -1 } + \alpha \,(\, 2 - a ^{1-\la } - a ^{\la -1 })
 - a- a\inv \ \Big] 
$$
while 
$
H_1(\la )_{ij} = (a-a\inv)\inv   \ \Big[ \ a ^{1-\la } - a ^{\la -1 } + 
\alpha (a ^{\la } - a ^{-\la })\ \Big] \ .
$ \EOE
\een
\end{rem}
\begin{cor}\label{descomposicion}
Given $N \in \orbu{D}$, consider the matrix decomposition 
$$
\dersin{\Delta_\la}{N}=
\begin{pmatrix}
\sub{A}{1N}(\la ) & 0\\
\sub{A}{2N}(\la ) & I   
\end{pmatrix}  \barr{l} \qn\\ I- \qn \earr \ ,
$$ 
as in Eq. \eqref{derivada en dos por dos}. 
Then $\|\sub{A}{1N}(\la )\|\leq \sub{k}{D,\la}\,$, where 
$$
\sub{k}{D,\la}  = \max\left\{ \max_{|d_i|\neq |d_j|}
\frac{|d_j|^{1-\la}|d_i|^\la+|d_i|^{1-\la}|d_j|^\la}{|d_i|+|d_j|}  ,  
\max_{\theta_i\neq \theta_j }|\la e^{i(\theta_j-\theta_i)}+1-\la| \right\}<1.
$$
\end{cor}
\proof 
As in the proof of Corollary 4.1.6 of \cite{[APS]}, it holds that 
$\|\sub{A}{1N}(\la )\|=\|\sub{A}{1D}(\la )\|$ for every  $N\in\orbu{D}$. 
On the other hand, by Proposition \ref{uno} and its notations, we get that 
$$
\displaystyle\|\sub{A}{1D}(\la )\|\leq \|\had{H_1(\la )}\|=\max_{ d_i\neq d_j}|H_1(\la )_{ij}| \ .
$$
If $|d_i| \neq |d_j|$, denote $a = \frac{\displaystyle|d_j|}{\displaystyle|d_i|} \neq 1$. 
As in Proposition \ref{uno}, we denote 
$$ 
A_{ij} =  |d_j|^{2-\la}|d_i|^\la-|d_i|^{2-\la}|d_j|^\la =
|d_j|\,|d_i| \ (a^{1-\la} - a^{\la-1}) \peso{and}  
$$
$$
B_{ij} = \ |d_j|^{1+\la}|d_i|^{1-\la}-|d_i|^{1+\la}|d_j|^{1-\la} = 
|d_j|\,|d_i| \ (a^{\la} - a^{-\la} )\ .
$$ 
Observe that $A_{ij}$ and $B_{ij}$  
have the same sign. So $|A_{ij} + B_{ij}| = \pm (A_{ij} + B_{ij}) $ and 
\begin{align*}
|H_1(\la )_{ij}| =  \left|\frac{ A_{ij} - e^{i(\theta_j-\theta_i)} B_{ij}}{|d_j|^2-|d_i|^2} \right| 
 &\leq \left|\frac{A_{ij} \ +\ B_{ij}}{|d_j|^2-|d_i|^2} \right| \\
 &=\frac{|d_j|-|d_i|}{|d_j|^2-|d_i|^2}\ \Big(|d_j|^{1-\la}|d_i|^\la+|d_i|^{1-\la}|d_j|^\la\Big)\\
&= \frac{|d_j|^{1-\la}|d_i|^\la+|d_i|^{1-\la}|d_j|^\la}{|d_j|+|d_i|} \ \ .
\end{align*}
This quantity is strictly lower that one (as observed in \cite{[HT]}) because
$$
\big(|d_j|+|d_i|\big)-\big(\big|d_j|^{1-\la}|d_i|^\la+|d_i|^{1-\la}|d_j|^\la\big)=
\big(|d_j|^\la-|d_i|^\la\big)\big(|d_j|^{1-\la}-|d_i|^{1-\la}\big)>0.
$$
On the other hand, if $|d_i|=|d_j|$ but $\theta_i\neq \theta_j\, $, using the 
triangle inequality we obtain that
$$
|H_1(\la )_{ij}|=\left|\la\big(e^{i(\theta_j-\theta_i)}-1\big)+1\right|=
\left|\la\,e^{i(\theta_j-\theta_i)}+(1-\la)\right|<1.
$$
In consequence, the bound for $\|\sub{A}{1N}(\la )\|$ is proved.
\QED

\subsection{The proof}

Now, we shall restate and prove Theorem \ref{EL TEO}:

\medskip
\noi
{\bf Theorem. } 
The $\la$-Aluthge transform $\alul{\cdot}:\orb{D}\to\orb{D}$ is a $C^\infty$ map, and for 
every $N\in\orbu{D}$, there exists a subspace $\sub{\cE}{N,\la}^s$ in the tangent space 
$\sub{T}{N}\orb{D}$ such that
\begin{enumerate}
	\item $\sub{T}{N}\orb{D}=\sub{\cE}{N,\la}^s\oplus \sub{T}{N}\orbu{D}$;
	\item Both, $\sub{\cE}{N,\la}^s$ and  $\sub{T}{N}\orbu{D}$, are $T_N\,\Delta_\la$-invariant;
	\item $\big\| T_N\,\Delta_\la \big| _{\sub{\cE}{N,\la}^s}\big\|\leq\sub{k}{D,\la}<1$, where 
	$\displaystyle \sub{k}{D,\la}$ is the constant of Corollary \ref{descomposicion}.
  \item If $U\in\matu$ satisfies $N=UDU^*$, then $\sub{\cE}{N,\la}^s=U(\sub{\cE}{D,\la}^s) U^*$.
  \end{enumerate}
  In particular, the map $\orbu{D} \ni N \mapsto \sub{\cE}{N,\la}^s$ is smooth.  
  This fact can be formulated in terms of the projections $\sub{P}{N,\la} $ onto 
  $\sub{\cE}{N,\la}^s$ parallel to $\sub{T}{N}\orbu{D}$, $N \in \orbu{D}$. \QED

\medskip
\noi
The proof of this theorem follows exactly the same steps as the proof of Theorem 3.1.1 of \cite{[APS]}, 
but now using Corollary \ref{descomposicion}, item \ref{pruyucciun}  of Remark \ref{remark sobre hadamard}, 
and Proposition \ref{uno} of this work. For this reason, we shall only give a sketch of the proof.

\proof[Sketch of proof] 
Fix $N=UDU^*\in\orbu{D}$. By the statement and the notations of 
Corollary \ref{descomposicion},  $\|\sub{A}{1N}(\la )\|<1$. So the operator 
$I-\sub{A}{1N}(\la )$ acting on $R(\qn)$ is invertible. Let $\cE^s_{N,\la}$ be the subspace defined by
$$
\cE^s_{N,\la}=\left\{\begin{pmatrix}
y\\
-\sub{A}{2N}(\la )(I-\sub{A}{1N}(\la )\,)^{-1} y
\end{pmatrix}:\ y\in R(\qn) \right\} \ ,
$$
Now, following the same steps as in the proof of Theorem 3.1.1 of \cite{[APS]} 
(with minor changes in order to adapts it to  our case) we can see that the following properties hold: 
\begin{itemize} 
\item The operator $\sub{P}{N,\la} \in L( \sub{T}{N}\orb{D} \, ) $ given by the matrix 
\beq\label{pnl}
\sub{P}{N,\la}=\begin{pmatrix}
I & 0\\
-\sub{A}{2N}(\la )(I-\sub{A}{1N}(\la )\,)^{-1} & 0
\end{pmatrix}
\begin{array}{cc}
\qn\\
I-\qn
\end{array}\
\end{equation}
is the projection onto $\cE_{N,\la}^s$ parallel to $\sub{T}{N}\orbu{D}$. 
Therefore we have the identity $\sub{T}{N}\orbu{D}=\cE^s_{N,\la} \oplus \sub{T}{N}\orbu{D} \,$. 
\item Since $T_N\,\Delta_\la=\sub{Ad}{U}(T_D\,\Delta_\la)\sub{Ad}{U}^{-1} $, then 
\beq\label{pnl2}
\sub{P}{N,\la}=\sub{Ad}{U}(\sub{P}{D,\la})\sub{Ad}{U}^{-1}  
\peso{and} \sub{P}{N,\la} ( \dersin{\Delta_\la}{N} ) = (\dersin{\Delta_\la}{N} ) \sub{P}{N,\la} \ .
\end{equation}

\item $\sub{\cE}{N,\la}^s=U(\sub{\cE}{D,\la}^s) U^*$ 
and both, $\cE^s_{N,\la}$ and $\sub{T}{N}\orbu{D}$, are invariant for 
$\dersin{\Delta_\la}{N}\,$.

\item $\big\| \big(\dersin{\Delta_\la}{N}\big)\big|_{\cE^s_{N,\la}} \big\| = 
\big\| \big(\dersin{\Delta_\la}{D}\big)\big|_{\cE^s_{D,\la}} \big\|$.  
So it suffices to show item 3 for $N = D$.

\item Let $Y=\begin{pmatrix} y\\ -\sub{A}{2D}(\la )(I-\sub{A}{1D}(\la )\,)^{-1} y \end{pmatrix} \in \cE^s_{D,\la}\,$, 
for some $y \in R(Q_D)$. Then
\begin{align*}
\|(\dersin{\Delta_\la}{D})\,Y\|_{_2}^2&=
\left\|\sub{A}{1D}(\la )y\right\|_{_2}^2 +\left\|\sub{A}{2D}(\la )y-\sub{A}{2D}(\la )
(I-\sub{A}{1D}(\la )\,)^{-1} y\right\|_{_2}^2\\&\leq \sub{k}{D,\la}^2
\left\|y\right\|_{_2}^2 +\left\|-\sub{A}{2D}(\la )\sub{A}{1D}(\la )
(I-\sub{A}{1D}(\la )\,)^{-1}y\right\|_{_2}^2 \ .
\end{align*}
where the inequality holds because $\|\sub{A}{1D}(\la )\|\leq\sub{k}{D,\la}\,$, 
by Corollary \ref{descomposicion}.

\item By item \ref{pruyucciun}  of Remark \ref{remark sobre hadamard} 
and Proposition \ref{uno},  we obtain that 
\begin{align*}
\left\|-\sub{A}{2D}(\la )\sub{A}{1D}(\la )(I-\sub{A}{1D}(\la )\,)^{-1} y\right\|_{_2}^2 & \le 
\sub{k}{D,\la}^2\left\|-\sub{A}{2D}(\la )(I-\sub{A}{1D}(\la )\,)^{-1} y\right\|_{_2}^2 \ .
\end{align*}

\item 
$\|(\dersin{\Delta_\la}{D})\,Y\|_{_2}^2 \leq \sub{k}{D,\la}^2\ \|y\|_{_2}^2+\sub{k}{D,\la}^2\ \left\|-
\sub{A}{2D}(\la )(I-\sub{A}{1D}(\la )\,)^{-1} y\right\|_{_2}^2=\sub{k}{D,\la}^2\|Y\|_{_2}^2\ .
$ Therefore we have that 
$\left.\big(\dersin{\Delta_\la}{N}\big)\right|_{\cE^s_{N,\la}}$ 
has norm lower or equal to $\sub{k}{D,\la}\,$. 
\end{itemize}
The smoothness of the map $\orbu{D} \ni N \mapsto \sub{\cE}{N,\la}^s$ 
follows from the existence of $C^\infty$ local cross sections  
for the map $\pi_D : \matu \to \orbu{D}$, stated  in Proposition \ref{son variedades}. 
For example, if $\sigma_D : \cV \to \matu$ is such a section near $D$, then by 
Eq. \eqref{pnl2}, 
\beq\label{pnl3}
\sub{P}{N,\la} = \sub{Ad}{\sigma_D (N)\,} \sub{P}{D,\la} \sub{Ad}{\sigma_D (N)^* \,} \quad , \quad N \in \cV \ .
\end{equation}
This completes the proof.
\QED

\begin{rem}\label{varia} \rm 
Using the notations and the proof of the Theorem, one can see that 
the significative parts of the projections $\sub{P}{N,\la}$ are 
\begin{align*}
(I- \qn )\sub{P}{N,\la} \qn & = -\sub{A}{2N}(\la )(I-\sub{A}{1N}(\la )\,)^{-1}  \\ 
& = - (I- \qn ) \dersin{\Delta_\la}{N} \qn  
\big(\qn - \qn \dersin{\Delta_\la}{N}  \qn \big)\inv  \ .
\end{align*}
When $N = D$, by  \ref{pruyucciun}  of Remark \ref{remark sobre hadamard} 
and  Proposition \ref{uno}, and using that the matrices 
$\{ 1 - H_1(\la )_{ij}\}_{ij} $ and  $\{ \big(1 - H_1(\la )_{ij} \big)\inv \}_{ij} $
are selfadjoint, we have that 
\begin{align}\label {uuu}
(I- \qd )\sub{P}{D,\la} \qd 
& = -(I- \qd ) \had{H_2(\la )} \qd  \big(\qd - \qd \had{H_1(\la )}  \qd \big)\inv  \nonumber \\
& = -(I- \qd ) \had{H_2(\la )}  \big( I - \had{H_1(\la )}\big)\inv  \qd 
\\& =  (I- \qd ) \had{G(\la )} \ \qd \ , \nonumber
\end{align}
where $G(\la )\in \matah$ has entries $G(\la )_{ij} = -H_2(\la )_{ij} \ (1 - H_1(\la )_{ij} )\inv $. 
Now, using  Remark \ref{con a}, we have the following properties, which have been 
 announced in Proposition \ref{no varian}:
\ben
\item If $|d_i | = |d_j|$, then  $G(\la )_{ij} = 0$. 
\item Suppose that all the eigenvalues of $D$ 
have the same moduli. Then, using Eq. \eqref{pnl} and Eq. \eqref{pnl2}, 
we get that  $\sub{P}{N,\la} = \qn $  (i.e., 
$\sub{\cE}{N,\la}^s = \big[ \sub{T}{N}\orbu{D} \big]^\bot$) for every $\la \in (0,1)$
and every $N \in \orbu{D}$.
\item If $|d_i | \neq |d_j|$,  $a = \frac{\displaystyle|d_j|}{\displaystyle|d_i|} \neq 1$ 
and $\beta = e^{i(\theta_j-\theta_i)} $, then
\beq \label{pnl5}
-G_{ij}(\la ) = \frac{ \ds a ^{1-\la } + a ^{\la -1 } +  
\beta \,(\, 2 - a ^{1-\la } - a ^{\la -1 }) - a- a\inv }
{\ds a-a\inv - ( a ^{1-\la } - a ^{\la -1 }) - \beta (a ^{\la } - a ^{-\la })}   = 
\frac {-z(\la)}{\ b(\la)} \ \ .
\end{equation}
Observe that $b(\la) \xrightarrow[\la \rightarrow 0]{} 0$, 
while $z(\la) \xrightarrow[\la \rightarrow 0]{} \beta (a + a\inv  -2)  \neq 0$. 
Therefore, we have that  $G(\la )_{ij} \xrightarrow[\la \rightarrow 0]{} \infty$. 
\item Suppose that  $D$ has at least two eigenvalues $d_i$ and $d_j$
such that $|d_i | \neq |d_j|$. 
Then, by the description $Q_D = \had {J} \pim \ \had {J} \inv $ given in \ref{pruyucciun}
 of Remark \ref{remark sobre hadamard}, and the
fact that  $G(\la )\in \matah$ (which implies that $\preal  \had{G(\la )}= \had{G(\la )}  \pim$), 
we have that 
$$
(I- \qd ) \had{G(\la )} \ \qd 
= \had {J} \preal  \had{G(\la )}  \pim \had {J} \inv =
\had {J}  \had{G(\la )}  \pim \had {J} \inv   \ .
$$ 
Since $G\in \matah \,$, then   $\|\had{G(\la )}  \| =\|\had{G(\la )}  \pim \|$. Hence 
$$
\|\sub{P}{D,\la} \| \ge \| (I- \qd )\sub{P}{D,\la} \qd \| = 
\| \had{G(\la )}  \| \ge |G(\la )_{ij}| \xrightarrow[\la \rightarrow 0]{} \infty  \ .
$$
Then the map $\la \mapsto \sub{P}{D,\la} $ can not be constant.  
\EOE
\een
\end{rem} 

\begin{rem}\label{pnl6}\rm
Note that, using Eq. \eqref{pnl3}, Remark \ref{varia} (particularly 
Eqs.  \eqref{uuu} and  \eqref{pnl5}\,), and the notations of the Theorem, 
we can conclude that the map 
$$
(0, 1) \times \orbu{D} \ni (\la , N ) \longmapsto  
\sub{Ad}{\sigma_D (N)\,} \sub{P}{D,\la} \sub{Ad}{\sigma_D (N)^* \,} = \sub{P}{N,\la} 
$$
is of class $C^\infty$. Another way to prove it is using that 
the map $(\la , N ) \mapsto \dersin{\Delta_\la}{N}$ is smooth, and 
then to apply Eqs. \eqref{derivada en dos por dos} and \eqref{pnl}. 
\EOE 
\end{rem}

\fontsize {9}{10}\selectfont

\end{document}